\documentclass[10pt,reqno]{amsart}
\usepackage{euscript,verbatim}
\usepackage{graphicx}%
\usepackage{multirow}%
\usepackage{amsmath,amssymb,amsfonts}%
\usepackage{amsthm}%
\usepackage{hyperref}

\usepackage{xcolor}
\newtheorem{theorem}{Theorem}[section]
\newtheorem*{theorem*}{Theorem}
\newtheorem{lemma}[theorem]{Lemma}

\theoremstyle{definition}
\newtheorem{definition}[theorem]{Definition}

\begin{document} 

\title[Gradient Bounds for Elliptic Equations]{Gradient Bounds for Viscosity Solutions to Certain Elliptic Equations} 

\author[T. Jeffres]{Thalia D. Jeffres} 
\address{Department of Mathematics and Statistics, Wichita State University, Wichita, KS, 67260-0033, USA}
\email{thalia.jeffres@wichita.edu}

\author[X. Li]{Xiaolong Li}\thanks{The second author's research is partially supported by NSF-DMS \#2405257 and a start-up grant at Auburn University}
\address{Department of Mathematics and Statistics, Auburn University, Auburn, AL, 36849, USA}
\email{xil0005@auburn.edu}




\begin{abstract}
Our principal object of study is the modulus of continuity of a periodic or 
uniformly vanishing function \( u: \mathbb{R} ^{n} \rightarrow \mathbb{R} \) which satisfies a 
degenerate 
elliptic equation \( F(x, u, \nabla u, D^{2} u) = 0 \) in the viscosity sense. The equations under consideration here have second-order terms of the form \( -{\rm Trace} \, (\mathcal{A} (\|\nabla u \|) \cdot D^{2} u) , \) where \( \mathcal{A} \) is an \( n\times n\) matrix which is symmetric and positive semi-definite. Following earlier work, \cite{Li21}, of the 
second author, which addressed the parabolic case, we identify a one-dimensional equation for which the modulus of continuity is a subsolution. In favorable cases, this one-dimensional operator can be used to derive a gradient bound on $u$ or to draw other conclusions about the nature of the solution. 
\end{abstract}

\maketitle


\section{Introduction} 

In this brief note we pursue a line of inquiry that was opened by the second author in \cite{Li21}. In that earlier work, the second author established that the modulus of continuity of a viscosity solution to certain fully nonlinear parabolic equations is a subsolution (also in the viscosity sense) of a corresponding one-dimensional equation. This one-dimensional operator arises from a structure condition first formulated by the second author \cite{Li21}. The structure condition is in turn inspired by an observation contained  in equations (3.11) through (3.15) and the surrounding discussion in \cite{CIL92}. Here we reformulate the structure condition in order to correct a minor error that occurred in \cite{Li21}. With this adjustment, the main result in \cite{Li21} remains valid. We also prove this result for solutions which vanish uniformly at infinity, and supply a few examples applying these results to the derivation of gradient bounds, when a gradient exists.  The methods and results of \cite{Li21} apply to parabolic equations; here we focus on the elliptic case. 

Our main result is the following, with definitions and concepts to be fully defined 
below. 

\begin{theorem}\label{thm main}
Suppose \( u: \mathbb{R} ^{n} \rightarrow \mathbb{R} \) is either a periodic or a 
uniformly vanishing viscosity solution to the degenerate elliptic equation \[ F(x,u, \nabla u, D^{2} u) = 0, \] where \( F: \mathbb{R}^n \times \mathbb{R} \times \mathbb{R} ^{n} \times \mathcal{S} (n) \rightarrow \mathbb{R} \) is a continuous function and \( \mathcal{S} (n) \) is the vector space of symmetric \( n\times n \) matrices.
Suppose that the pair \( (F,f) \) satisfies the structure condition in Definition \ref{def SC}. Then the modulus of
continuity, \( \omega ,\) of $u$ is a viscosity subsolution to the equation \[ f(s, \phi , \phi ', \phi '') =0 \] on \( (0, \infty ). \)   
\end{theorem}

If the one-dimensional operator satisfies a Comparison Principle, and if a supersolution of \( f =0 \) having suitable shape can be found, then further implications such as gradient bounds follow. Statements and examples are given in Section 5. 

The utility of the modulus of continuity lies in its applicability to functions lacking differentiability. For this reason, it is well adapted to the viscosity setting and under favorable conditions can even serve as a substitute for the gradient. Such connections were explored by Andrews and Xiong, \cite{AX19}, who obtained gradient estimates for a 
large class of quasilinear elliptic  equations. Their results precede ours and are 
also more general. On the other hand, their method is different and part of our purpose here is to illustrate a computational technique, showing how to use the structural condition to identify a one-dimensional operator and then to use that to draw conclusions about the nature of the solution to the original equation. 

In a series of papers, Andrews and Clutterbuck made an 
ingenious application of the modulus of continuity, using  it to prove the Fundamental Gap Conjecture, \cite{AC11},  and to determine a sharp lower bound on the first eigenvalue of the Sch\"{o}dinger operator, \cite{AC13}. 
The connection is through heat equations and separation of variables. Similar methods were used again by these two authors, \cite{AC09a, AC09}, to establish bounds on the spatial gradient of a solution to a parabolic equation. Later, Ni \cite{Ni13} gave a proof of the fundamental gap estimate that remained entirely within the elliptic setting. Modulus of continuity is a powerful tool which continues to find new applications. In \cite{LeBalch21}, Le Balc'h makes applications to special cases of the Landis Conjecture. For more related works, see \cite{AN12}, \cite{Andrewssurvey15}, \cite{Li16, Li21}, \cite{SWW19}, \cite{AX19}, \cite{HWZ20}, \cite{DSW21}, \cite{LW17, LW21, LW21JGA, LW21MRL, LW23}, \cite{LTW24}, \cite{WZ23, WZ24}, and the references therein.

Going back further, all of the works mentioned above use a doubling method originated by Kruzhkov, \cite{Kruzkov67}, and a two-point maximum principle for functions of the form 
\( u(x) - u(y) - \phi (x,y). \)  A standard approach in proving gradient bounds for quasilinear parabolic equations is first to prove them over the boundary, and then, via maximum principle, in the interior. Kruzhkov discovered that in some cases, the gradient could be bounded in one step, that in doubling the number of spatial variables one introduces a kind of internal boundary along the diagonal. A gradient bound over this boundary is then a gradient bound over the original region. (In addition to the original reference, see also the text of Lieberman, \cite{Lieberman96book}, for discussion of this method.) In fact, Andrews and Xiong, \cite{AX19}, work directly with this two-point maximum principle rather than with the modulus of continuity itself. 

What is very interesting is that functions of the form 
\( u(x) - u(y) - \phi (x,y) \) 
are also extremely important in the development of viscosity solutions, as in \cite{CIL92} and also 
summarized briefly in Section 2 of the present paper. 
In the context of viscosity solutions, these functions arise for reasons that are very different from Kruzhkov's motivations. A viscosity solution is continuous, and therefore attains a maximum over any compact set, but the sub- and superjets used to define a viscosity solution could at that point contain nothing helpful or even be empty. These are the objects that stand in for first and second derivatives, and without them, analysis of the shape of the graph - in other words, maximum and comparison principles - is impossible. The existence result, Theorem 3.2 of \cite{CIL92}, shows us that the relevant superjets at local maxima of 
\( u(x) -u(y) - \phi (x,y) \) contain elements with properties like those of literal gradients and Hessians. 
Perhaps it is this connection that explains the deep compatibility that emerges among 
modulus of continuity, two-point maximum principles, and viscosity solutions.

This article is organized as follows. In Section 2, we recall the definitions of viscosity solutions and recall the maximum principle for semi-continuous functions in \cite{CIL92}. In Section 3, we formulate the structure condition and present the proof of Theorem \ref{thm main} for periodic solutions. The case of uniformly vanishing solutions is then proved in Section 4. In Section 5, we discuss applications of the modulus of continuity estimates to gradient bounds with two examples and conclude the paper with a general strategy summarized in Theorem \ref{thm gradient}. 

\section{Preliminaries.} 

We provide in this section the relevant definitions for viscosity solutions, and also record the important 
existence result of Crandall, Ishii, and Lions, \cite{CIL92}. That foundational paper gives a full 
exposition of the concepts, vocabulary, and results essential to viscosity solutions. 

The equations under consideration here have the form 
\[ F(x, u,  \nabla u, D^{2} u) = 0, \] 
where $F$ is a function 
\[ F: \Omega \times \mathbb{R} \times \mathbb{R} ^{n} \times \mathcal{S} (n) \rightarrow \mathbb{R} , \] 
where  \( \Omega \subset \mathbb{R} ^{n} \) is an open subset of \( \mathbb{R} ^{n} \) and \( \mathcal{S} (n) \) is the vector space of symmetric \( n\times n \) matrices. $F$ is assumed 
to be continuous, and the typical element of the domain is written \( (x, z, \vec{p} , A). \) $F$ is called degenerate elliptic if it is order-reversing in the matrix position, meaning that if  \( A \leq B, \) then 
\[ F(x, z,\vec{p} ,B ) \leq F(x, z, \vec{p} ,A). \] 

Suppose that $u$ is upper semi-continuous in a subset \( \mathcal{O} \) of \( \mathbb{R} ^{n} \) 
which contains a point \( x_{0} ,\) and that \( \phi \) is a \( C^{2} \) function in a neighborhood of \( x_{0} .\)  If \( u-\phi \) attains a local maximum at \( x_{0} , \) then for \( x\in \mathcal{O} \) and near \( x_{0} , \) 
\[ u(x) \leq u(x_{0} ) + \langle \nabla \phi (x_{0} ) , x-x_{0} \rangle + \frac{1}{2} \langle D^{2} \phi (x_{0} ) \cdot (x-x_{0} ) , x-x_{0} \rangle + r(x; x_{0} ) , \] 
with 
\[ \lim _{x\rightarrow x_{0} } \frac{r(x;x_{0} ) }{\| x-x_{0} \| ^{2} } =0. \] 
Taking this as motivation, the pair \( (\vec{p} , A) \in \mathbb{R} ^{n} \times \mathcal{S} (n) \) belongs to the second-order superjet of $u$ at \( x_{0} \) if for $x$ near \( x_{0} \) we  have 
\[ u(x) \leq u(x_{0} ) + \langle \vec{p} , x-x_{0} \rangle + \frac{1}{2} \langle A \cdot (x-x_{0} ) , x -x_{0} \rangle + r(x;x_{0} ), \] 
with 
\[ \lim _{x\rightarrow x_{0} , x\in \mathcal{O} } \frac{r(x;x_{0} )}{\| x-x_{0} \| ^{2} } =0.  \] 
The second-order superjet of $u$ at \( x_{0} , \) denoted \( (J^{2, +} _{\mathcal{O} } u)(x_{0} ) , \) is the union of all such. 

\begin{definition}
    An upper semicontinuous function \( f: \mathcal{O} \rightarrow \mathbb{R} \) is a 
{\em viscosity subsolution} of \( F =0 \) if \( F(x, u(x), \vec{p} ,A ) \leq 0 \) for all 
\( x\in \mathcal{O} \) and \( (\vec{p} , A) \in (J^{2, +} _{\mathcal{O} } u)(x). \) 
\end{definition}


The corresponding notions of subjets and supersolutions are obtained by reversing the direction of the inequalities. Namely, a pair \( (\vec{p} , A) \in \mathbb{R} ^{n} \times \mathcal{S} (n) \) lies in the second-order subjet of $u$ at \( x_{0} \) if for $x$ near \( x_{0}  \in \mathcal{O} , \) 
the inequality 
\[ u(x) \geq u(x_{0} ) + \langle \vec{p} , x-x_{0} \rangle + \frac{1}{2} \langle A \cdot (x-x_{0} ) , x -x_{0} \rangle + r(x;x_{0} ), \] 
with 
\[ \lim _{x\rightarrow x_{0} , x\in \mathcal{O} } \frac{r(x;x_{0} )}{\| x-x_{0} \| ^{2} } =0,   \] 
holds, and the second-order subjet \( (J^{2,-} _{\mathcal{O} } u) (x_{0} ) \) is the union of all such. 

\begin{definition}
    A lower semicontinuous function \( f: \mathcal{O} \rightarrow \mathbb{R} \) is a 
{\em viscosity supersolution } of \( F =0 \) if \( F(x, u(x), \vec{p} ,A ) \geq 0 \) for all 
\( x\in \mathcal{O} \) and \( (\vec{p} , A) \in (J^{2, -} _{\mathcal{O} } u)(x). \) 
\end{definition}

A continuous function \( f: \mathcal{O} \rightarrow \mathbb{R}  \) is a viscosity solution of 
\( F =0 \) if it is both a viscosity subsolution and supersolution. In particular, by being 
both upper- and lower semicontinuous, $u$ will be  continuous. 

For our purposes here, \( \mathcal{O} \) will always be an open subset of \( \mathbb{R} ^{n} . \) 
In this case, \( (J^{2,+} _{\mathcal{O} } u) (x_{0} ) \) for \( x_{0} \in \mathcal{O} \) is the set of pairs \( (\nabla \phi (x_{0} ) , D^{2} \phi (x_{0} ) ) \) for which \( \phi \in C^{2} (U) \) for 
\( U \subset \mathcal{O} \) an open set containing \( x_{0} , \) and \( u-\phi \) has a local 
maximum at \( x_{0} . \) When \( \mathcal{O} \) is open, one usually writes just \( (J^{2, + } u) (x_{0} ) . \) In the case in which \( \mathcal{O} \) is open, a third characterization is possible, and we will use this in Sections 3 and 4. One can in this case say that \( (J^{2, + } _{\mathcal{O} } u )(x_{0} ) \) is the set of pairs \( (\nabla \phi (x_{0} ), D^{2} \phi (x_{0} ) ) \) for which  \( \phi \in C^{2} (U) \) for 
\( U \subset \mathcal{O} \) an open set containing \( x_{0}  , \) and \( \phi (x_{0} ) = u(x_{0} ), \) with \( \phi (x) \geq u(x) \) throughout $U.$ 

We include now the statement of the existence result of \cite{CIL92}, adapted to our purposes. This is Theorem 3.2 of \cite{CIL92}, where it is stated in fuller generality. 

\begin{theorem}[Theorem (3.2 of \cite{CIL92}]\label{thm CIL}
Let $u$ be a continuous function in an open set \( \mathcal{O} \subset \mathbb{R} ^{n}, \) and \( \phi \in C^{2} (U) \) be defined in an open set \( U \subset \mathbb{R} ^{n} \times \mathbb{R} ^{n} \) containing \( \mathcal{O} \times \mathcal{O} . \) Suppose that 
\( (x_{0} , y_{0} ) \in \mathcal{O} \times \mathcal{O} \) is a local maximum of the function of 
two variables \( u(x) - u(y) - \phi (x,y). \) Let \( \varepsilon > 0 \) be any positive real number. Corresponding to \( \varepsilon , \) there exist symmetric matrices \( A, B \in \mathcal{S} (n) \) 
having the following properties. 
\begin{enumerate}
    \item The pair \( (\nabla _{x} \phi (x_{0} , y_{0} ), A ) \) belongs to the closure of the superjet of $u$ at \( x_{0} , \, (\overline{J} ^{2,+} u) (x_{0} ) . \) 

\item The pair \( (\nabla _{y} \phi (x_{0} , y_{0} ) , -B) \) belongs to \( (\overline{J} ^{2, +} (-u)) (y_{0} ). \) 

\item The matrices $A$ and $B$ satisfy 
\[ -\left( \frac{1}{\varepsilon } + \| D^{2} \phi (x_{0} , y_{0} ) \| \right)  \cdot I \leq \left[ \begin{array}{rr} 
  A & 0 \\ 
  0 & -B 
 \end{array} \right] \leq D^{2} \phi (x_{0} , y_{0} ) + \varepsilon (D^{2} \phi (x_{0} ,y_{0} ) )^{2} . \] 
\end{enumerate}

\end{theorem}

 
This can be seen as a sort of first- and second-derivative test, as was discussed previously. 

To a viscosity solution  $u$ of the equation \( F =0 \) we associate a function of a single real 
variable, its modulus of continuity. 

\begin{definition}
    Let \( u: \Omega \rightarrow \mathbb{R} \) be a given function on \( \Omega \subset \mathbb{R} ^{n} . \) If $h$ satisfies 
\[ | u(x) -u(y) | \leq 2 h(\frac{\| x-y \| }{2} ) \] 
for all \( x, y \in \Omega , \) then $h$ is called a {\em modulus of continuity.} The (optimal) {\em modulus of continuity} \( \omega \) of $u$ is given by 
\[ \omega (s) = \sup \left\{ \frac{1}{2} (u(x) -u(y) ) | x, y\in \Omega , \, \| x-y \| = 2s \right\} .\] 
\end{definition}
 

\section{Formulation of the structure condition and proof of the main theorem} 

Many naturally occurring differential operators contain second-order terms of the form 
\[ -{\rm Trace} (\mathcal{A} (x, u, \nabla u) \cdot D^{2} u), \] 
with \( \mathcal{A} \) an \( n\times n \) matrix that is symmetric and positive semidefinite. For the moment, let us examine  just on this second-order part, because this will motivate the 
structure condition. For this second-order part, the  
corresponding function 
\[ F: \Omega \times \mathbb{R} \times \mathbb{R} ^{n} \times \mathcal{S} (n) \rightarrow \mathbb{R} \] 
is then \( F(x,z, \vec{p} ,A) = -{\rm Trace} \, (\mathcal{A} (x,z, \vec{p} ) \cdot A). \) This important class of operators includes  the minimal surface operator, the $p-$Laplace operator, and 
the Laplacian itself. For these operators, \(  \mathcal{A} (x, z,\vec{p}) \geq 0, \) and this 
implies the degenerate ellipticity condition that 
\[ F(x, z,\vec{p} ,B) \leq F(x,z, \vec{p} ,A) \] 
whenever \( A \leq B. \) 

Asking what is true when $A$ and $B$ are those matrices whose existence is asserted in the existence result of \cite{CIL92}, Theorem 3.2,  leads us to a structure condition. Part of the answer to this question is provided in \cite{Li21}, where this result appeared. 
\begin{lemma}[Lemma 4.1 of \cite{Li21}]
    If  \( A, B \in \mathcal{S} (n) \) satisfy 
\[ \left[ \begin{array}{rr} 
    A & 0 \\ 
    0 & -B 
    \end{array} \right] \leq D^{2} (2\phi (\frac{\| x-y\| }{2} ) ) , \] 
then \( A \leq B, \) and 
\[ {\rm Trace} \, (A-B) \leq 2 \phi ''(\frac{\| x-y\| }{2} ). \] 
\end{lemma}

What is true if $A$ and $B$ satisfy the inequality that appears on the right hand side of the Existence Result, Theorem 3.2, of \cite{CIL92}? Namely,  if  for a given \( \varepsilon > 0, \) the inequality 
\[ \left[ \begin{array}{rr} 
  A & 0 \\ 
  0 & -B 
 \end{array} \right] \leq D^{2} ( 2\phi (\frac{\| x-y \| }{2} ) )  + \varepsilon \left( D^{2} ( 2\phi (\frac{\| x-y \| )}{2} ) ) \right)  ^{2}  \] 
 holds, what is true? 
(Incidentally, the factors of $2$ and $1/2$ appear because the specific way in which the \( C^{2} \) function of \cite{CIL92} depends on $x$ and $y$ in \( \mathbb{R} ^{n} \) is through \( 2\phi (\| x-y \| /2), \) where \( \phi \) is a function of a single real variable.) 
Consistently with the notation of \cite{Li21}, we write 
\[ P = 2 D_{x} ^{2} \phi (\frac{\| x-y\| }{2} ) , \] 
and then find that if \( Q = P + 2\varepsilon P^{2} , \) then 
\[ \left[ \begin{array}{rr} 
      A & 0 \\ 
      0 & -B 
     \end{array} \right] \leq \left[ \begin{array}{rr} 
                      Q & -Q \\ 
                      -Q & Q 
                      \end{array} \right] . \] 
From this it follows that \( A -B \leq 0. \) 

Following the same reasoning as in Li's proof of lemma 4.1 of \cite{Li21}, and also of \cite{AC09}, where  
a similar calculation is done, we note that if $C$ is any symmetric \( n\times n \) matrix for which 
\[ \left[ \begin{array}{cc} 
    I & C \\ 
    C & I 
    \end{array} \right] > 0, \] 
then \( {\rm Trace} (A-B) \leq 2 \cdot {\rm Trace} (I-C)Q. \) Making the specific choice 
\( C = I - 2 u \cdot u^{T} \) for \( \vec{u} \) a unit vector, a calculation then yields 
\[ {\rm Trace } (A-B) \leq 2 \phi ''(\frac{\| x-y\| }{2} ) + 8 \varepsilon \| P\vec{u} \| ^{2} . \] 
We summarize these findings. 
\begin{lemma}
    If \( A, B \in \mathcal{S} (n) \) satisfy 
\[ \left[ \begin{array}{rr} 
  A & 0 \\ 
  0 & -B 
 \end{array} \right] \leq D^{2} ( 2\phi (\frac{\| x-y \| }{2} ) )  + \varepsilon \left( D^{2} ( 2\phi (\frac{\| x-y \| )}{2} ) ) \right)  ^{2} , \] 
then \( A \leq B, \) and 
\[ {\rm Trace } (A-B) \leq 2 \phi ''(\frac{\| x-y\| }{2} ) + 8 \varepsilon \| P\vec{u} \| ^{2} . \] 
\end{lemma}

For the Laplace operator, \( \mathcal{A} = I, \) and we will then have \( F(x,z,\vec{p} ,A) = -{\rm Trace} A, \) so that 
\[ F(y, z, \vec{p} ,B) - F(x, w,\vec{p} ,A) = {\rm Trace} (A-B) \leq 2 \phi ''(\frac{\| x-y\| }{2} ) + 8 \varepsilon \| P\vec{u} \| ^{2} . \]
For the minimal surface operator, the $p$-Laplacian, or other operators of the form \( -{\rm Trace} (\mathcal{A} (\| \vec{p} \| ) \cdot A) \) with \( \mathcal{A} \geq 0, \) we have 
\[ F(y,z,\vec{p} ,B) - F(x, w, \vec{p} ,A) \leq \lambda (\| \vec{p} \| ) \left[ 2\phi '' (\frac{\| x-y\| }{2} ) + 8\varepsilon (\| P\vec{u} \| ^{2} \right]  , \] 
with \( \lambda (\| \vec{p} \| ) \) denoting the minimal eigenvalue of \( \mathcal{A} (\| \vec{p} \| ) . \) Examination of these and other examples motivates the formulation of the structure condition. 
This is slightly different from the structure condition that was given in \cite{Li21}, because of the inclusion of the quadratic terms. We now allow $F$ to again be more general, as described in Section 2. 

\begin{definition}[\bf The structure condition]\label{def SC}
    Let \( \varepsilon > 0\) be given, along with \( x, y \in \mathbb{R} ^{n} \) with \( \| x-y\| = 2s >0, \) and \( z,w \in \mathbb{R} \) with \( w-z = 2\phi (s) \geq 0, \) and \( A, B \in \mathcal{S} (n) \) satisfying 
\[ \left[ \begin{array}{rr} 
  A & 0 \\ 
  0 & -B 
 \end{array} \right] \leq D^{2} ( 2\phi (\frac{\| x-y \| }{2} ) )  + \varepsilon \left( D^{2} ( 2\phi (\frac{\| x-y \| )}{2} ) ) \right)  ^{2} . \] 
 Examine \( F(y,z,\vec{p} ,B) - F(x,w, \vec{p} ,A) \) for the choice \( \vec{p} = \phi ' \cdot \vec{u} , \) where 
\[ \vec{u} = \frac{x-y}{\| x-y\| } . \] 
If 
\[ F(y,z,\vec{p} ,B) - F(x,w,\vec{p} ,A) \leq -2 f(s,\phi ,\phi ' , \phi '') + \varepsilon q(s,\phi ,\phi ', \phi '') , \]
for a function $q$ depending only on the indicated quantities, 
then the pair \( (F,f) \) satisfies the structure condition. 
\end{definition}

What is worth emphasizing here is that the additional term arising from the quadratic term in the existence result, the term containing the factor \( \varepsilon ,\) does not depend on $A$ and $B.$ This is important because in the existence result of \cite{CIL92}, the matrices $A$ and $B$ themselves 
depend on \( \varepsilon ,\) so that varying \( \varepsilon \) causes $A$ and $B$ to change. This independence of $q$ from $A$ and $B$ means that the main theorem of \cite{Li21}, Theorem 1.1, and its proof remain true with only small adjustments, as we will now see. We will address the uniformly vanishing case separately, 
in Section 4. 

\begin{theorem}\label{thm 3.4}
Suppose that \( u: \mathbb{R} ^{n} \rightarrow \mathbb{R} \) is a periodic viscosity solution to 
\[ F(x,u, \nabla u, D^{2} u) =0, \] 
where 
\( F: \mathbb{R}^n \times \mathbb{R} \times \mathbb{R} ^{n} \times \mathcal{S} (n) \rightarrow \mathbb{R} \) is continuous and degenerate elliptic in the sense defined in Section 2. 
Suppose that the pair \( (F,f) \) satisfies the structure condition in Definition \ref{def SC}. Then \( \omega ,\) the modulus of
continuity  of $u$ is a viscosity subsolution to the equation \[ f(s, \phi , \phi ', \phi '') =0 \] on \( (0, \infty ). \) 
\end{theorem}


Periodic solutions to elliptic equations arise naturally in many ways. In \cite{Post03}, for example, Post 
considers the noncompact case, in which a spectral gap need not occur at all. This author 
uses periodicity to construct manifolds with as many spectral gaps as desired.

\begin{proof}[Proof of Theorem 3.4] Let \( s_{0} \in (0, \infty ). \) Since \( (0, \infty ) \) is an open set, we need to show that if \( \phi \) is any \( C^{2} \) function defined in a neighborhood \( U \subset (0, \infty ) \) which contains \( s_{0} , \) with \( \omega (s_{0} ) = \phi (s_{0} ) \) and 
satisfying \( \phi \geq \omega ,\) then 
\[ f(s_{0} , \phi (s_{0} ), \phi '(s_{0} ) , \phi ''(s_{0} ) ) \leq 0. \] 
Let \( \delta > 0 \) be small enough that \( (s_{0} -\delta , s_{0} + \delta ) \) is contained in the domain of \( \phi . \) Then if points $x$ and $y$ of \( \mathbb{R} ^{n} \) satisfy 
\( \| x-y\| = 2s, \) for \( s\in (s_{0} -\delta , s_{0} + \delta ) , \) the very definition of \( \omega \) together with the properties describing \( \phi \) imply that the function of two variables 
\[ u(x) -u(y) - 2 \phi (\frac{\| x-y\| }{2} ) \] 
is non-positive, 
\[ u(x) -u(y) - 2 \phi (\frac{\| x-y\| }{2} ) \leq 0. \] 
Moreover, the periodicity and continuity of $u$ imply that zero is an attained maximum value of this function, that there are points \( x_{0}  , y_{0} \in \mathbb{R} ^{n} \) with \( \| x_{0} -y_{0} \| = 2s_{0} \) and 
\[ u(x_{0} ) -u(y_{0} ) - 2 \phi (\frac{\| x_{0} -y_{0} \| }{2} ) = 0. \] 
The existence result, Theorem 3.2 of \cite{CIL92} or Theorem \ref{thm CIL}, now asserts that for a given choice of \( \varepsilon > 0,\) there exist 
corresponding  matrices \( A, B \in \mathcal{S} (n) \) such that with 
\[ \vec{u} = \frac{x_{0} -y_{0} }{\| x_{0} -y_{0} \| } , \] 
that 
\[ (\phi '(s_{0} ) \vec{u} , A) \in (\overline{J} ^{2,+} u)(x_{0} ) , \] 
and 
\[ (\phi '(s_{0} ) \vec{u} , -B) \in (\overline{J} ^{2,+} (-u))(y_{0} ) , \] 
and 
\[ \left[ \begin{array}{rr} 
  A & 0 \\ 
  0 & -B 
 \end{array} \right] \leq D^{2} ( 2\phi (\frac{\| x_{0} -y_{0}  \| }{2} ) )  + \varepsilon \left( D^{2} ( 2\phi (\frac{\| x_{0} -y_{0}  \| )}{2} ) ) \right)  ^{2} . \] 
Using this and what it means to be a viscosity solution, it follows that 
\[ F(x_{0}, u(x_{0} ), \phi '(s_{0} ) \vec{u} , A) - F(y_{0} ,u(y_{0} ), \phi '(s_{0} ) \vec{u} , B) \leq 0. \] 
Together with the structure condition in Definition \ref{def SC}, we conclude that 
\begin{align*}
     0 & \leq F(y_{0} , u(y_{0} ) , \phi '(s_{0} ) \vec{u} ,B) - F(x_{0} , u(x_{0} ) , \phi '(s_{0} ), \vec{u} , A) \\
     & \leq -2f(s_{0} , \phi , \phi ' , \phi '') + \varepsilon q(s_{0} , \phi , \phi ' , \phi '') , 
\end{align*}
with \( \phi \) and its derivatives evaluated at \( s_{0} .\) 
Varying \( \varepsilon \) can change $A$ and $B$ in the intermediate term, but through suppressing 
this  term we have simply 
\[ 0  \leq -2f(s_{0} , \phi , \phi ' , \phi '') + \varepsilon q(s_{0} , \phi , \phi ' , \phi '') , \] 
and $A$ and $B$ do not appear. Letting \( \varepsilon \rightarrow 0, \) it must be that 
\[ f(s_{0} ,  \phi (s_{0} ), \phi '(s_{0} ), \phi ''(s_{0} ) )\leq 0, \] 
which is what we wanted to show.
\end{proof}

In Section 5, we will demonstrate how to identify the one-dimensional operator $f$ and in favorable cases how to use the knowledge that \( \omega \) is a subsolution of \( f=0 \) to 
draw conclusions.

\section{Proof of the main theorem in the case of a uniformly vanishing solution} 

Suppose that \( u: \mathbb{R} ^{n} \rightarrow \mathbb{R} \) is a solution in the viscosity 
sense to the equation 
\[ F(x, u, \nabla u, D^{2} u) =0, \] 
for $F$ as described in Section 2, and that $u$ is uniformly vanishing at infinity, meaning that for any \( \varepsilon >0, \) there 
exists a corresponding \( R > 0  \) for which \( | u(x) | < \varepsilon \) whenever \( \| x \| >R. \) In this situation, too, if the pair \( (F,f) \) satisfies the structure condition, then the modulus of continuity of $u$ is a subsolution to \( f=0 \) on \( (0, \infty ). \) 

Suppose \( s_{0} \in (0, \infty ), \) and that \( \phi \) is a \( C^{2} \) function satisfying 
\( \phi (s_{0} ) = \omega (s_{0} ) \) and \( \omega (s) \leq \phi (s) \) throughout a neighborhood 
\( (s_{0} - \delta , s_{0} + \delta ), \) with \( \delta >0 \) chosen so that \( (s_{0} - \delta , s_{0} + \delta ) \subset (0, \infty ). \)  The function of two variables 
\[ u(x) -u(y) - 2 \phi (\frac{\| x-y\| }{2} ) \] 
is defined for all points  \( (x,y) \in \mathbb{R} ^{n} \times \mathbb{R} ^{n} \) for which 
\[ \frac{\| x-y\| }{2} \in (s_{0} -\delta , s_{0} + \delta ). \] 
This set is open and unbounded, but the function nevertheless has an achieved maximum value of zero, for the following reason.  Restricted to the closed set 
\( \{ \| x-y\| = 2s_{0} \} , \, u(x) - u(y) \) is a bounded function, and \( u(x) -u(y) \rightarrow 0 \) if both \( \| x\| \rightarrow \infty \) and \( \| y\| \rightarrow \infty . \) This  means that among points for which \( \| x-y\| = 2s_{0} ,  \) the function \( u(x) - u(y)  \) has an achieved maximum.  So there is a point \( (x_{0} ,y_{0} ) \) with \( \| x_{0} -y_{0} \| = 2s_{0} , \) and 
\[ u(x_{0} ) - u(y_{0} ) = 2\omega (s_{0} ). \] 
Since \( \phi (s_{0} ) = \omega (s_{0} ), \) it follows that 
\[ u(x_{0} ) - u(y_{0} ) - 2\phi (\frac{\| x_{0} -y_{0} \| }{2} ) =0. \] 
At other points  \(  (x,y) \in \mathbb{R} ^{n} \times \mathbb{R} ^{n} \) for which \(  \| x-y\| = 2s, \) with \(  s\in (s_{0} -\delta , s_{0} + \delta )  , \) we have 
\[ u(x) -u(y) \leq 2\omega (s) \leq 2 \phi (s)  . \] 
This  means that at \( (x_{0} ,y_{0} ), \) the function of two variables has an achieved maximum 
value of zero. The existence theorem 3.2 of \cite{CIL92} can now be applied and the remainder of the argument 
is the same as in Section 3 above. We record this result. 
\begin{theorem}
Suppose that \( F: \mathbb{R}^n \times \mathbb{R} \times \mathbb{R} ^{n} \times \mathcal{S} (n) \rightarrow \mathbb{R} , \)  is a continuous function 
which is degenerate elliptic 
and that \( u: \mathbb{R} ^{n} \rightarrow \mathbb{R} \) is a viscosity solution to 
\[ F(x, u, \nabla u, D^{2} u) =0 \] 
that vanishes uniformly at infinity. 
Suppose that the pair \( (F,f) \) satisfies the structure condition in Definition \ref{def SC}. Then \( \omega ,\) the modulus of 
continuity  of $u$ is a viscosity subsolution to the equation \[ f(s, \phi  , \phi ', \phi '') =0 \] on \( (0, \infty ). \) 
\end{theorem}


\section{Applications to Gradient Bounds.} 

If $u$ is a harmonic function, then the norm of its gradient  and geometrically meaningful combinations such as \( \sqrt{1 + \| \nabla u \| ^{2} }  \)  
are subharmonic; similar statements hold for solutions to the heat equation. This simple calculation serves as the prototype for the Bernstein method. One identifies a closely related equation for which an expression in \( \| \nabla u \| \) is a subsolution, and then applies the Comparison Principle. The paper of Serrin, \cite{Serrin71}, is a survey of this technique in the elliptic and parabolic settings.

When solutions are merely continuous and therefore needn't possess literal gradients, the modulus of continuity might serve  as a substitute. This is a function of one real variable only. Employing the Bernstein strategy, the insight of Li, \cite{Li21}, is to identify a one-dimensional operator for which this function of a single real variable is a subsolution. Here in the elliptic case, we  have seen in Sections 3 and 4 that when the pair \( (F,f) \) satisfies the structure condition, \( f(s, \omega , \omega ', \omega '' ) \leq 0.\) If a supersolution \( \zeta \) satisfying \( \omega \leq \zeta \) on the boundary of the one-dimensional domain can be found, and if $f$ satisfies a Comparison Principle strong enough to apply in the viscosity setting, then \( \omega \leq \zeta \) throughout the entire domain. We illustrate how to use this knowledge  with two examples. 

\vspace{.5cm} 

{\bf Example 1.} Suppose \( F: \mathbb{R} ^{n} \times \mathbb{R} \times \mathbb{R} ^{n} \times \mathcal{S} (n) \rightarrow \mathbb{R} \) is given by 
\[ F(x, z,\vec{p} ,A) = -{\rm Trace} (\mathcal{A} (\vec{p} ) \cdot A) + \langle \vec{b} (x), \vec{p} \rangle + cz, \] 
where \( \vec{b} \) is a bounded vector field, $c$ is a positive constant, and \( \mathcal{A} (\vec{p} ) \) is a symmetric matrix for which these exists \( \lambda > 0 \) with the property that 
\( \mathcal{A} (\vec{p} ) \geq \lambda I \) holds for all \( \vec{p} \in \mathbb{R} ^{n} . \) 

If \( u: \mathbb{R} ^{n} \rightarrow \mathbb{R} \) is either a periodic solution or a uniformly 
vanishing solution and if \( \omega \) is its modulus of continuity, then \( f(\omega ) \leq 0 \) 
for the choice 
\[ f(\phi ) = -\lambda \phi '' -B | \phi '| + c\phi , \] 
where $B$ is an upper bound on the norm of the vector field, \( \| \vec{b} (x) \| \leq B. \) Non-negative solutions to the equation \( f= 0\) satisfying \( \zeta (0) =0 \) take the form 
\[ \zeta (s) = \mu ^{2} (e^{\alpha _{1} s } - e^{\alpha _{2} s} ) , \] 
where 
\[ \alpha _{1} = \frac{ -B +\sqrt{B^{2} + 4\lambda c} }{2\lambda } , \] 
which is always positive, and 
\[ \alpha _{2} = \frac{-B -\sqrt{B^{2} + 4\lambda c} }{2\lambda } , \] 
which is always negative. 
For the periodic case, if $D$ is the diameter of a compact repeating region, then \( \omega (D/2) \leq \| u\| _{0} , \) and so if \( \mu ^{2} \) is chosen so that 
\[ \mu ^{2} (e^{\alpha _{1} D/2 } -e^{\alpha _{2} D/2} ) \geq \| u\| _{0} , \] 
then \( \omega (D/2) \leq \zeta (D/2). \) 

In the case of a uniformly vanishing solution, \( \omega (s) \rightarrow \| u\| _{0} /2 \) as \( s\rightarrow \infty , \) so if any \( \varepsilon > 0 \) is chosen, there corresponds a value \( s(\varepsilon ) >0 \) with the property that 
\[ \omega (s) \leq \frac{\| u\| _{0} }{2} + \varepsilon \] 
for all \( s> s(\varepsilon ). \) In this case, then, choosing \( \mu ^{2} \) large enough that 
\[ \mu ^{2} (e^{\alpha _{1} s(\varepsilon ) } -e^{\alpha _{2} s(\varepsilon ) } ) \geq \frac{\| u\| _{0} }{2} + \varepsilon \] 
will ensure that \( \omega (s(\varepsilon ) ) \leq \zeta (s(\varepsilon ) ). \) 

Applying the Comparison Principle, (3.14) from \cite{CIL92}, we conclude that \( \omega (s) \leq \zeta (s) \) throughout the entire interval \( [0, D/2] \) or \( [0, s(\varepsilon ) ] \) of \( \mathbb{R}\).

Replacing the original repeating region by a larger one, say of diameter \( k\cdot D/2 \)  for 
\( k\in \mathbb{Z} ^{+} , \) the constant \( \mu ^{2} \) can now be chosen smaller, because the increasing and unbounded function \( \zeta \) has to clear the height \( \| u\| _{0} \) at the larger value \( s = k\cdot D/2,  \) 
which is to say that \( \mu ^{2}  \) must satisfy only 
\[ \mu ^{2} (e^{\alpha _{1} k D/2 } - e^{\alpha _{2} k D/2} ) \geq \| u\| _{0} .\] 
On the original interval, \( (0, D/2), \) the function 
\( \omega \) must lie below all of them, even as \( \mu ^{2} \rightarrow 0, \) and so \( \omega (s) \equiv 0. \) 

The same argument applies to the case of a uniformly vanishing solution.  Taking the interval to be larger, say \( [0, s(\varepsilon ) + a] \) for \( a>0, \) allows 
\( \mu ^{2} \) to be chosen smaller, in fact arbitrarily so, as \( a\rightarrow \infty ,  \) 
because \( \mu ^{2} \) must now satisfy 
\[ \mu ^{2} (e^{\alpha _{1} (s(\varepsilon ) +a)} - e^{\alpha _{2} (s(\varepsilon ) +a) } ) \geq \frac{\| u\| _{0} }{2} + \varepsilon .\] 
By the 
Comparison Principle, \( \omega \) lies below all of them on the fixed interval \( [0, s(\varepsilon ) ] ,\) and this is only possible if \( \omega \equiv 0 \) on \( [0, s(\varepsilon )] . \) 

In either case, $u$ must be constant in balls of radius \( D/2,  \) for the case of a periodic 
solution, or \( s(\varepsilon ), \) for the case of a uniformly vanishing solution, and hence everywhere. From the original equation, it follows that $u$ is identically zero. This is consistent with direct application of the Maximum Principle in case of a \( C^{2} \) solution.

It isn't necessary that the equation \( f =0 \) be solved exactly, as was possible in the previous 
example. A supersolution of suitable shape suffices. 

{\bf Example 2.} Suppose \( F(x,z, \vec{p} ,A) = - {\rm Trace} (\mathcal{A} (\| \vec{p} \|  ) \cdot A) + g(\vec{p} ) + cz, \) where \( \mathcal{A} (\| \vec{p} \|  ) \) is symmetric and positive-semidefinite and $c$ is a positive constant. 

This is a large class of operators which includes those with the minimal surface operator or the $p$-Laplace operator as second-order terms.

Then we identify the one-dimensional operator 
\[ f(\phi ) = - \lambda (|\phi '| ) \phi '' + c\phi . \] 
If \( \mathcal{A} (\| \vec{p} \| )  \) arises from the minimal surface operator, then 
\[ f(\phi ) = - \frac{\phi ''}{(1+ |\phi '| )^{3/2} } + c\phi , \] 
and for the $p$-Laplace operator, 
\[ f(\phi ) = -| \phi '| ^{p-2} \phi '' + c\phi . \] 
If $u$ is a periodic solution with a compact repeating region of diameter $D,$ then if $a$ is 
chosen so that \( a\geq \| u\| _{0} , \) then the downward-opening parabola with 
axis of symmetry at \( D/2 , \) 
\[ \zeta (s) = -\frac{4a}{D^{2} } (s-\frac{D}{2} )^{2} + a \] 
satisfies the differential inequality \( f(\zeta ) \geq 0 \) with 
\[ \limsup  _{\| x-y\|  \rightarrow 0}  \frac{u(x) -u(y) }{\| x-y\| } \leq \frac{4a}{D} . \] 
By taking larger repeating regions, we conclude that a periodic solution must be zero. By a similar argument, the same conclusion follows if $u$ is a uniformly vanishing solution. 

\vspace{.5cm} 

Extracting from these examples some general principles, if \( \omega (x) \leq \zeta (s) \) can be shown to hold on the domain of \( \omega ,\) then drawing further conclusions depends on the nature of the comparison function \( \zeta .\) The proof of the following statement is the same as in 
the examples. 
\begin{theorem}\label{thm gradient}
Suppose that \( u: \mathbb{R} ^{n} \rightarrow \mathbb{R} \) is either a periodic or a uniformly vanishing solution to the degenerate elliptic equation 
\[ F(x, u, \nabla u, D^{2} u) =0 \] 
in the viscosity sense, with $F$ a continuous function 
\[ F: \mathbb{R} ^{n} \times \mathbb{R} \times \mathbb{R} ^{n} \times \mathcal{S} (n) \rightarrow \mathbb{R}. \] 
If $f$ is a one-dimensional operator for which the pair \( (F,f) \) satisfies the structure 
condition and the Comparison Principle (3.14) of \cite{CIL92}, and if \( \zeta \) is a supersolution 
of \( f=0 \) with \( \omega \leq \zeta \) on the boundary of the domain of \( \omega , \) then 
the following is true. 

{\bf (A)} If \( \zeta \) is bounded, then for all \( x, y \in \mathbb{R} ^{n} \) with \( \| x-y \| \leq D \) (if $u$ is periodic), or \( \| x-y \| \leq 2s(\varepsilon ) \) (if $u$ is uniformly vanishing), 
\[ |u(x) - u(y)| \leq 2 \cdot \sup \zeta. \] 

{\bf (B)} If \( \zeta \) is differentiable at zero, with \( \zeta (0) =0, \) then $u$ is Lipschitz continuous, and also 
\[ \limsup _{\| x-y\| \rightarrow 0}  \frac{ | u(x) -u(y) | }{\| x-y \| } \leq \zeta '(0). \] 
If at some particular  \( x\in \mathbb{R} ^{n} \) the gradient \( \nabla u \) exists, then 
\[ \lim _{y \rightarrow x} \frac{| u(x) - u(y) | }{\| x-y \| } \leq \zeta '(0), \] 
and this implies a bound on \( \| \nabla u \|  \) at that point. 

{\bf (C)} If for some \( \alpha \in (0, 1], \, \zeta \) has an expansion of the form 
\[ \zeta (s) = 0 + a s^{\alpha } + r(s;0) \] 
with 
\[ \lim _{s\rightarrow 0} \frac{r(s;0) }{s^{\alpha } } = 0, \] 
then $u$ is uniformly H\"{o}lder continuous. 
\end{theorem}

In scenario (A), it is always true that \( |u(x) -u(y)|  \leq 2 \| u\| _{0} , \) and 
\( \zeta \) depends on \( \| u\| _{0} , \) but this does not mean that the statement is vacuous. The function \( \zeta \) must be larger than \( \| u\| _{0} ,\) for periodic solutions, or than \( \| u \| _{0} /2 + \varepsilon \) for uniformly 
vanishing solutions,  only at the right hand endpoint of \( [0, D/2] \) or \( [0, s(\varepsilon ) ] \) respectively. Moreover, the examples illustrate how the choice of domain and of \( \zeta \) are not unique, and that making different choices might allow one to draw further conclusions. In the situations described in (B) and (C), since the construction of \( \zeta \) depends on \( \| u \| _{0} , \) one obtains a gradient or H\"{o}lder quotient bound which depends on \( \| u \| _{0} , \) and it is very typical first to obtain estimates on the solution itself, and then on the gradient in terms of the \( C^{0} \) bound.

\section{Acknowledgments} Each author held a postdoctoral position at the University of California at Irvine. Though separated by many years, we each learned a great deal from Professor Peter Li and from the many in his orbit, most especially Professor Lei Ni. It is a privilege to acknowledge here their erudition and generosity.

\section{Declarations} 

{\bf Conflict of Interest.} Not applicable (no potential conflict of interest). 

\section{Data Availability Statement}
Data sharing is not applicable to this article as no datasets were generated or analyzed during the current study.

\bibliographystyle{plain}
\bibliography{myref}

\begin{thebibliography}{10}

\bibitem{Andrewssurvey15}
Ben Andrews.
\newblock Moduli of continuity, isoperimetric profiles, and multi-point estimates in geometric heat equations.
\newblock In {\em Surveys in differential geometry 2014. {R}egularity and evolution of nonlinear equations}, volume~19 of {\em Surv. Differ. Geom.}, pages 1--47. Int. Press, Somerville, MA, 2015.

\bibitem{AC09a}
Ben Andrews and Julie Clutterbuck.
\newblock Lipschitz bounds for solutions of quasilinear parabolic equations in one space variable.
\newblock {\em J. Differential Equations}, 246(11):4268--4283, 2009.

\bibitem{AC09}
Ben Andrews and Julie Clutterbuck.
\newblock Time-interior gradient estimates for quasilinear parabolic equations.
\newblock {\em Indiana Univ. Math. J.}, 58(1):351--380, 2009.

\bibitem{AC11}
Ben Andrews and Julie Clutterbuck.
\newblock Proof of the fundamental gap conjecture.
\newblock {\em J. Amer. Math. Soc.}, 24(3):899--916, 2011.

\bibitem{AC13}
Ben Andrews and Julie Clutterbuck.
\newblock Sharp modulus of continuity for parabolic equations on manifolds and lower bounds for the first eigenvalue.
\newblock {\em Anal. PDE}, 6(5):1013--1024, 2013.

\bibitem{AN12}
Ben Andrews and Lei Ni.
\newblock Eigenvalue comparison on {B}akry-{E}mery manifolds.
\newblock {\em Comm. Partial Differential Equations}, 37(11):2081--2092, 2012.

\bibitem{AX19}
Ben Andrews and Changwei Xiong.
\newblock Gradient estimates via two-point functions for elliptic equations on manifolds.
\newblock {\em Adv. Math.}, 349:1151--1197, 2019.

\bibitem{CIL92}
Michael~G. Crandall, Hitoshi Ishii, and Pierre-Louis Lions.
\newblock User's guide to viscosity solutions of second order partial differential equations.
\newblock {\em Bull. Amer. Math. Soc. (N.S.)}, 27(1):1--67, 1992.

\bibitem{DSW21}
Xianzhe Dai, Shoo Seto, and Guofang Wei.
\newblock Fundamental gap estimate for convex domains on sphere---the case {$n = 2$}.
\newblock {\em Comm. Anal. Geom.}, 29(5):1095--1125, 2021.

\bibitem{HWZ20}
Chenxu He, Guofang Wei, and Qi~S. Zhang.
\newblock Fundamental gap of convex domains in the spheres.
\newblock {\em Amer. J. Math.}, 142(4):1161--1191, 2020.

\bibitem{Kruzkov67}
S.~N. Kru\v{z}kov.
\newblock Nonlinear parabolic equations with two independent variables.
\newblock {\em Trudy Moskov. Mat. Ob\v{s}\v{c}.}, 16:329--346, 1967.

\bibitem{LeBalch21}
K\'evin Le~Balc'h.
\newblock Exponential bounds for gradient of solutions to linear elliptic and parabolic equations.
\newblock {\em J. Funct. Anal.}, 281(5):Paper No. 109094, 29, 2021.

\bibitem{Li16}
Xiaolong Li.
\newblock Moduli of continuity for viscosity solutions.
\newblock {\em Proc. Amer. Math. Soc.}, 144(4):1717--1724, 2016.

\bibitem{Li21}
Xiaolong Li.
\newblock Modulus of continuity estimates for fully nonlinear parabolic equations.
\newblock {\em Calc. Var. Partial Differential Equations}, 60(5):Paper No. 182, 23, 2021.

\bibitem{LTW24}
Xiaolong Li, Yucheng Tu, and Kui Wang.
\newblock On a class of quasilinear operators on smooth metric measure spaces.
\newblock {\em Comm. Anal. Geom.}, 32(10):2757--2803, 2024.

\bibitem{LW17}
Xiaolong Li and Kui Wang.
\newblock Moduli of continuity for viscosity solutions on manifolds.
\newblock {\em J. Geom. Anal.}, 27(1):557--576, 2017.

\bibitem{LW21}
Xiaolong Li and Kui Wang.
\newblock Lower bounds for the first eigenvalue of the {L}aplacian on {K}\"ahler manifolds.
\newblock {\em Trans. Amer. Math. Soc.}, 374(11):8081--8099, 2021.

\bibitem{LW21JGA}
Xiaolong Li and Kui Wang.
\newblock Sharp lower bound for the first eigenvalue of the weighted {$p$}-{L}aplacian {I}.
\newblock {\em J. Geom. Anal.}, 31(8):8686--8708, 2021.

\bibitem{LW21MRL}
Xiaolong Li and Kui Wang.
\newblock Sharp lower bound for the first eigenvalue of the weighted {$p$}-{L}aplacian {II}.
\newblock {\em Math. Res. Lett.}, 28(5):1459--1479, 2021.

\bibitem{LW23}
Xiaolong Li and Kui Wang.
\newblock Eigenvalue estimates on quaternion-{K}\"ahler manifolds.
\newblock {\em J. Geom. Anal.}, 33(3):Paper No. 85, 20, 2023.

\bibitem{Lieberman96book}
Gary~M. Lieberman.
\newblock {\em Second order parabolic differential equations}.
\newblock World Scientific Publishing Co., Inc., River Edge, NJ, 1996.

\bibitem{Ni13}
Lei Ni.
\newblock Estimates on the modulus of expansion for vector fields solving nonlinear equations.
\newblock {\em J. Math. Pures Appl. (9)}, 99(1):1--16, 2013.

\bibitem{Post03}
Olaf Post.
\newblock Periodic manifolds with spectral gaps.
\newblock {\em J. Differential Equations}, 187(1):23--45, 2003.

\bibitem{Serrin71}
James Serrin.
\newblock Gradient estimates for solutions of nonlinear elliptic and parabolic equations.
\newblock In {\em Contributions to nonlinear functional analysis ({P}roc. {S}ympos., {M}ath. {R}es. {C}enter, {U}niv. {W}isconsin, {M}adison, {W}is., 1971)}, pages 565--601. Academic Press, New York-London, 1971.

\bibitem{SWW19}
Shoo Seto, Lili Wang, and Guofang Wei.
\newblock Sharp fundamental gap estimate on convex domains of sphere.
\newblock {\em J. Differential Geom.}, 112(2):347--389, 2019.

\bibitem{WZ23}
Kui Wang and Shaoheng Zhang.
\newblock Lower bounds for the first eigenvalue of {$p$}-{L}aplacian on {K}\"ahler manifolds.
\newblock {\em Proc. Amer. Math. Soc.}, 151(6):2503--2515, 2023.

\bibitem{WZ24}
Kui Wang and Shaoheng Zhang.
\newblock Lower bounds for the first eigenvalue of the {$p$}-{L}aplacian on quaternionic {K}\"ahler manifolds.
\newblock {\em J. Geom. Anal.}, 34(5):Paper No. 123, 19, 2024.

\end{thebibliography}
\end{document}